\magnification= \magstep1
\input amstex
\documentstyle{amsppt}
\define\defeq{\overset{\text{def}}\to=}

\define\Gal{\operatorname{Gal}}

\def \isom {\overset \sim \to \rightarrow}

\define\Spec{\operatorname{Spec}}

\def \Ker{\operatorname {Ker}}
\def \id{\operatorname {id}}
\def \Fr{\operatorname {Fr}}
\def \char{\operatorname {char}}
\def \Spec{\operatorname {Spec}}
\def \Spf{\operatorname {Spf}}
\def \Sp{\operatorname {Sp}}
\def \geo{\operatorname {geo}}

\def \cd{\operatorname {cd}}

\def \nor{\operatorname {nor}}

\def \Ind{\operatorname {Ind}}
\def \rig{\operatorname {rig}}

\NoRunningHeads
\NoBlackBoxes
\topmatter

\title 
On \'etale fundamental groups of formal fibres of $p$-adic 
curves
\endtitle
\bigskip

\author
Mohamed Sa\"\i di
\bigskip
\endauthor

\abstract We investigate a certain class of (geometric) finite (Galois) coverings of formal fibres of $p$-adic curves and the corresponding quotient of the (geometric) \'etale fundamental group. 
A key result in our investigation is that these (Galois) coverings can be compactified to finite (Galois) coverings of proper $p$-adic curves. We also prove that the maximal prime-to-$p$ quotient of the geometric \'etale fundamental group of a (geometrically connected) formal fibre of a $p$-adic curve is (pro-)prime-to-$p$ {\it free of finite computable rank}. 
\endabstract

\toc

\subhead
\S0. Introduction/Main Results
\endsubhead

\subhead
\S1. Background
\endsubhead

\subhead
\S2. Geometric Galois groups of formal boundaries of formal germs of $p$-adic formal curves
\endsubhead

\subhead
\S3. Geometric fundamental groups of formal fibres of $p$-adic curves
\endsubhead

\endtoc

\endtopmatter

\footnotetext [1]{2010 Mathematics Subject Classification. Primary 14H30; Secondary 11G20.}
\footnotetext [2]{Key words and phrases. Formal fibre, $p$-adic curves, \'Etale fundamental groups.}

\document 

\subhead
\S 0. Introduction/Main Results
\endsubhead
A classical result in the theory of \'etale fundamental groups is the description of the structure of the geometric \'etale fundamental group of 
an affine, smooth, and geometrically connected curve over a field of characteristic $0$ (cf. [Grothendieck], Expos\'e XIII, Corollaire 2.12). 
In this paper we investigate the structure of a certain quotient of the geometric \'etale fundamental group of a formal fibre of a $p$-adic 
curve.

Let $R$ be a complete discrete valuation ring, $K=\Fr (R)$ its quotient field, and
$k$ its residue field which we assume to be algebraically closed of characteristic $p\ge 0$. 
Let $X$ be a proper, flat, and normal formal $R$-curve whose special fibre $X_k$ is reduced and consists of $n\ge 1$ distinct
irreducible components $\{P_{i}\}_{i=1}^n$ which intersect at a (closed) point $x\in X_k(k)$, and $x$ is the unique singular point of $X_k$. 
Write $\widetilde P_i\to P_i$ for the morphism of normalisation. We assume $\widetilde P_i=\Bbb P^1_k$ is a projective line, the morphism $\widetilde P_i\to P_i$ is a homeomorphism, and if $x_i$ is the (unique) pre-image of $x$ in $\widetilde P_i$ then $x_i\in \widetilde P_i(k)$ is the zero point of $\widetilde P_i$. 
In particular, the configuration of the irreducible components of $X_k$ is tree-like. 
The formal curve $X$ has a formal covering $X=B\cup D_1\cup \ldots \cup D_n$ where
$B\subset X$ is a formal sub-scheme with special fibre $B_k=X_k\setminus \{\infty_i\}_{i=1}^n$ ($\infty_i$ is the image in $P_i$ of the 
infinity point of $\widetilde P_i$, $1\le i\le n$),
$D_i=\Spf \langle \frac {1}{T_i}\rangle$ is an $R$-formal closed unit disc with special fibre
$D_{i,k}=P_i\setminus \{x\}$ and generic fibre $D_{i,K}=\Sp K\langle \frac {1}{T_i}\rangle$ which is a closed unit $K$-rigid disc centred 
at the point $\infty_i\in D_{i,K}(K)$ (which specialises in $\infty _i\in D_{i,k}$), $1\le i\le n$.
Write $\Cal F\defeq {\Cal F}_x=\Spf \hat {\Cal O}_{X,x}$
for the formal germ of $X$ at $x$ and 
$\Cal F_K\defeq \Cal F_{x,K}=\Spec (\hat \Cal O_{X,x}\otimes _RK)$ for the formal fibre of the generic fibre $X_K$ of 
the algebraisation of $X$ at $x$ 
(cf. 1.2 for more details, as well as Remark 3.1 which asserts that any formal germ of a 
formal $R$-curve at a closed point admits a compactification as above).

Let $S\subset \Cal F_K$ be a (possibly empty) finite set of closed points. Write $\pi_1(\Cal F_K\setminus S)^{\geo}$ for the
geometric \'etale fundamental group of $\Cal F_K\setminus S$ (in the sense of Grothendieck, cf. 1.3 for more details), 
and consider the quotient $\pi_1(\Cal F_K\setminus S)^{\geo}\twoheadrightarrow \widehat \pi_1(\Cal F_K\setminus S)^{\geo}$
which classifies finite coverings $\Cal Y'\to \Cal F'\defeq \Cal F\times_R{R'}$, where $R'/R$ is a finite extension, $\Cal Y'$ is normal and geometrically connected, 
which are \'etale above $\Cal F_{K'}\setminus S_{K'}$ ($K'\defeq \Fr R'$ and $S_{K'}\defeq S\times_KK'$) and \'etale above the generic points of $\Cal F_k$ (cf. loc. cit.).
Similarly, write $\pi_1(X_K\setminus (S\cup \{\infty_i\}_{i=1}^n))^{\geo}$ for the geometric \'etale fundamental group of the affine curve $X_K\setminus (S\cup \{\infty_i\}_{i=1}^n)$
and consider the quotient $\pi_1(X_K\setminus (S\cup \{\infty_i\}_{i=1}^n))^{\geo}\twoheadrightarrow \widehat \pi_1(X_K\setminus (S\cup \{\infty_i\}_{i=1}^n);\{\infty_i\}_{i=1}^n)^{\geo}$ 
which classifies finite coverings $Y'\to X_{R'}$ which are \'etale above $X_{K'}\setminus (S_{K'}\cup \{\infty_i\}_{i=1}^n)$, possibly ramified above the points 
$\{\infty_i\}_{i=1}^n$ with ramification indices prime-to-$p$, and which are \'etale above the generic points of $X_k$ (here $R'$, $K'$ and $S_{K'}$ are as above). We also write 
$\widehat \pi_1(X_K\setminus (S\cup \{\infty_i\}_{i=1}^n))^{\geo,p}\defeq \widehat \pi_1(X_K\setminus (S\cup \{\infty_i\}_{i=1}^n);\{\infty_i\}_{i=1}^n)^{\geo,p}$ 
(resp. $\widehat \pi_1(\Cal F_K\setminus S)^{\geo,p}$) for the maximal pro-$p$ quotient
of  $\widehat \pi_1(X_K\setminus (S\cup \{\infty_i\}_{i=1}^n);\{\infty_i\}_{i=1}^n)^{\geo}$
(resp. $\widehat \pi_1(\Cal F_K\setminus S)^{\geo}$).
Our first main result is the following (cf. Theorem 3.2).

\proclaim {Theorem 1} The (scheme) morphism $\Cal F_{K}\to X_K$ induces a continuous homomorphism $\widehat \pi_1(\Cal F_{K}\setminus S)^{\geo}
\to \widehat \pi_1(X_K\setminus (S\cup \{\infty_i\}_{i=1}^n);\{\infty_i\}_{i=1}^n)^{\geo}$ (resp. $\widehat \pi_1(\Cal F_{K}\setminus S)^{\geo,p}
\to \widehat \pi_1(X_K\setminus (S\cup \{\infty_i\}_{i=1}^n);\{\infty_i\}_{i=1}^n)^{\geo,p}$)
which makes $\widehat \pi_1(\Cal F_{K}\setminus S)^{\geo}$ (resp. $\widehat \pi_1(\Cal F_{K}\setminus S)^{\geo,p}$) into a semi-direct factor (cf. Definition 1.1.4 and Lemma 1.1.5) of 
$\widehat \pi_1(X_K\setminus (S\cup \{\infty_i\}_{i=1}^n); \{\infty_i\}_{i=1}^n)^{\geo}$ (resp. $\widehat \pi_1(X_K\setminus (S\cup \{\infty_i\}_{i=1}^n);\{\infty_i\}_{i=1}^n)^{\geo,p}$). 
In particular, the above homomorphisms are injective.
\endproclaim

In the course of proving Theorem 1 (cf. proof of Theorem 3.2) we prove the following.
\proclaim{Theorem 2} Let $f:\Cal Y\to \Cal F$ be a finite (Galois) covering with $\Cal Y$ normal and geometrically connected, 
which is \'etale above $\Cal F_K\setminus S$ and above the generic points of $\Cal F_k$. Then there exists, after possibly a finite extension of $K$, a finite (Galois) 
covering $\tilde f: Y\to X$ of formal schemes with $Y$ normal and geometrically connected, which is \'etale above $X_K\setminus (S\cup\{\infty_i\}_{i=1}^n)$ and above the generic points of $X_k$, 
is possibly ramified above the points $\{\infty_i\}_{i=1}^n$ with ramification indices prime-to-$p$, and which induces by pull back via the (scheme) morphism $\Cal F\to X$ the covering $f$.
\endproclaim


Let $g_x\defeq \text{genus}(X_K)$, which is also called the genus of the formal fibre $\Cal F_K$.
Write $\pi_1(\Cal F_K\setminus S,\eta)^{\geo,p'}$ (resp. $\pi_1(X_K\setminus (S\cup \{\infty_i\}_{i=1}^n),\eta)^{\geo,p'}$)
for the maximal prime-to-$p$ quotient of the geometric \'etale fundamental group $\pi_1(\Cal F_K\setminus S,\eta)^{\geo}$ (resp. $\pi_1(X_K\setminus (S\cup \{\infty_i\}_{i=1}^n),\eta)^{\geo}$).
Our second main result is the following (cf. Theorem 3.4).

\proclaim {Theorem 3} Let $S(\overline K)=\{y_1,\ldots,y_m\}$ of cardinality $m\ge 0$. Then the continuous homomorphism 
$\pi_1(\Cal F_K\setminus S,\eta)^{\geo,p'}\to \pi_1(X_K\setminus (S\cup \{\infty_i\}_{i=1}^n),\eta)^{\geo,p'}$ (induced by the (scheme) morphism $\Cal F_K\to X_K$) 
is an isomorphism. In particular,
$\pi_1(\Cal F_K\setminus S,\eta)^{\geo,p'}$ is (pro-)prime-to-$p$
free of rank $2g_x+n+m-1$ and can be generated by $2g_x+n+m$ generators 
$\{a_1,\ldots,a_g,b_1,\cdots,b_g,\sigma_1,\ldots,\sigma_n,\tau_1,\ldots,\tau_m\}$
subject to the unique relation $\prod_{i=1}^g[a_i,b_i]\prod _{j=1}^n\sigma _j\prod_{t=1}^m \tau _t=1$, where $\sigma_j$ (resp $\tau _t$) is a generator of inertia at $\infty _i$
(resp. $y_t$). 
\endproclaim

Next, we outline the content of the paper. In $\S1$ we collect some well-known background material.
In $\S2$ we investigate a certain quotient of the absolute Galois group of a formal boundary of a formal germ 
of a $p$-adic curve and prove Proposition 2.5 which is used in the proof of Theorem 1. In $\S3$ we prove Theorems 1 and 3.

\subhead
Notations
\endsubhead
In this paper $K$ is a complete discrete valuation field, $R$ its valuation ring, $\pi$ a uniformising parameter, 
and $k\defeq R/\pi R$ the residue field which we assume to be algebraically closed of characteristic $p\ge 0$. 

We refer to [Raynaud], $3$, for the terminology we will use 
concerning $K$-rigid analytic spaces, $R$-formal schemes, as well as the link between formal and rigid geometry. 
For an $R$-(formal) scheme $X$ we will denote by
$X_K\defeq X\times _RK$ (resp. $X_k\defeq X\times _Rk$) the generic (resp. special) fibre of $X$ 
(the generic fibre is understood in the rigid analytic sense in the case where $X$ is a formal scheme).
Moreover, if $X=\Spf A$ 
is an affine formal $R$-scheme of finite type we denote by $X_K\defeq \Sp (A\otimes _RK)$ the associated $K$-rigid affinoid space and
will also denote, when there is no risk of confusion, by $X_K$
the affine scheme $X_K\defeq \Spec (A\otimes _RK)$.

A formal (resp. algebraic) $R$-curve is an $R$-formal scheme of finite type (resp. $R$-scheme of finite type)
flat, separated, and whose special fibre is equi-dimensional of dimension $1$.
For a $K$-scheme (resp. $K$-rigid analytic space) $X$ and $L/K$ a field extension (resp. a finite extension) we write
$X_L\defeq X\times _KL$ which is an $L$-scheme (resp. an $L$-rigid analytic space). If $X$ is a proper and normal
formal $R$-curve we also denote, when there is no risk of confusion, by $X$ the algebraisation of $X$ which is an algebraic $R$-curve and by $X_K$ 
the proper normal and algebraic $K$-curve associated to the rigid $K$-curve $X_K$ via the rigid GAGA functor.

For a profinite group $H$ and a prime integer $\ell$ we denote by $H^{\ell}$ the maximal pro-$\ell$ quotient of $H$,
and $H^{\ell'}$ the maximal prime-to-$\ell$ quotient of $H$.


\subhead
\S 1 Background
\endsubhead
In this section we collect some background material used in this paper.

\subhead
1.1 
\endsubhead
Let $p>1$ be a prime integer. We recall some well-known facts 
on profinite pro-$p$ groups.
First, we recall the following characterisations of free pro-$p$ groups.

\proclaim {Proposition 1.1.1} Let $G$ be a profinite pro-$p$ group. Then the following properties are equivalent.

\noindent
(i)\ $G$ is a free pro-$p$ group.

\noindent
(ii)\ The $p$-cohomological dimension of $G$ satisfies $\cd_p(G)\le 1$.

\noindent
In particular, a closed subgroup of a free pro-$p$ group is free.

\endproclaim

\demo{Proof} Well-known (cf. [Serre], and [Ribes-Zalesskii], Theorem 7.7.4).
\qed
\enddemo

Next, we recall the notion of a {\it direct factor} of a free pro-$p$ group 
(cf. [Garuti], 1, the discussion preceding  Proposition 1.8, and [Sa\"\i di], $\S1$).

\proclaim {Definition/Lemma 1.1.2 (Direct factors of free  pro-$p$ groups)} 
Let $F$ be a free pro-$p$ group, $H\subseteq F$ a closed 
subgroup, and $\iota: H\to F$ the natural homomorphism. We say that 
$H$ is a direct factor of $F$ if there exists a continuous homomorphism 
$s:F\to H$ such that $s\circ \iota=\id_H$. 
There exists then a (non unique) closed subgroup $N$ of $F$ such that $F$ is isomorphic to the free direct product $H\star N$.
We will refer to such a subgroup $N$ as a supplement of $H$.
\endproclaim

\demo{Proof}  See [Sa\"\i di] Lemma 1.1.2.
\qed
\enddemo

One has the following cohomological characterisation of direct factors of free 
pro-$p$ groups.

\proclaim {Proposition 1.1.3}
Let $H$ be a pro-$p$ group, $F$ a free pro-$p$ group, and 
$\sigma :H\to F$ a continuous homomorphism. Assume that the map induced by $\sigma$ on 
cohomology
$$h^1(\sigma):H^1(F,\Bbb Z/p\Bbb Z)\to H^1(H,\Bbb Z/p\Bbb Z)$$
is surjective, where
$\Bbb Z/p\Bbb Z$ is considered as a trivial discrete module. 
Then $\sigma$ induces an isomorphism $H\isom \sigma(H)$ and $\sigma(H)$ is a direct factor of $F$. 
In particular, $H$ is pro-$p$ free.
We say that $\sigma$ makes $H$ into a direct factor of $F$.
\endproclaim

\demo{Proof} cf. [Garuti], Proposition 1.8 and Proposition 1.1.1 above.
\qed
\enddemo

Next, we consider the notion of a {\it semi-direct factor} of a profinite group.

\definition {Definition 1.1.4 (Semi-direct factors of profinite groups)} Let $G$ be a profinite group, $H\subseteq G$ a closed subgroup,
and $\iota: H\to G$ the natural homomorphism. We say that 
$H$ is a {\it semi-direct factor} of $G$ if there exists a continuous homomorphism 
$s:G\to H$ such that $s\circ \iota=\id_H$ ($s$ is necessarily surjective). 
\enddefinition

\proclaim{Lemma 1.1.5} Let $\tau:H\to G$ be a continuous homomorphism between profinite groups. Write 
$H=\underset{j\in J} \to {\varprojlim} H_j$ as the projective limit of the inverse system $\{H_j,\phi_{j'j},J\}$
of finite quotients $H_j$ of $H$ with index set $J$.
Suppose there exists,  $\forall j\in J$, a surjective homomorphism $\psi_j :G\twoheadrightarrow H_j$ such that 
$\tau \circ \psi_j:H\twoheadrightarrow H_j$ is the natural map and
$\psi_j=\phi _{j'j}\circ \psi _{j'}$ whenever this makes sense.
Then $\tau$ induces an isomorphism $H\isom \tau(H)$ and $\tau(H)$ is a semi-direct factor of $G$. 
We say that $\tau$ makes $H$ into a semi-direct factor of $G$.
\endproclaim

\demo{Proof} Indeed, the $\{\psi_j\}_{j\in J}$ give rise to a continuous (necessarily surjective) homomorphism $\psi :G\to H$ which is a right inverse of $\tau$.
\qed
\enddemo

\subhead
1.2. Formal Patching
\endsubhead
Next, we explain the procedure which allows to construct (Galois) coverings
of curves in the setting of formal geometry by patching coverings of formal (affine, non-proper)
curves with coverings of formal germs at closed points of the special fibre along the boundaries of these formal germs.

\subhead {1.2.1} 
\endsubhead
Let $X$ be a proper, normal, formal $R$-curve with $X_k$ reduced.
For $x\in X$ a closed point
let $\Cal F_x\defeq \Spf \hat {\Cal O}_{X,x}$ be the formal completion of $X$ at $x$ which we will refer to as the {\it formal germ} of $X$ at $x$.
Thus, $\hat {\Cal O}_{X,x}$ is the completion of the local ring of the algebraisation of $X$ at $x$.
Write $\Cal F_{x,K}\defeq \Spec (\hat {\Cal O}_{X,x}\otimes _RK)$.
We will refer to $\Cal F_{x,K}$ as the {\it formal fibre} of $X_K$ at $x$. 
Let $\{\Cal P_i\}_{i=1}^n$ be the minimal prime ideals of $\hat {\Cal O}_{X,x}$ which contain $\pi$; they correspond 
to the branches $\{\eta_i\}_{i=1}^n$ of the completion of $X_k$ at $x$ (i.e., closed points of the normalisation of $X_k$ above $x$), 
and $\Cal X_i=\Cal X_{x,i}\defeq \Spf \hat {\Cal O}_{x,\Cal P_i}$
the formal completion of the localisation of $\Cal F_x$ at $\Cal P_i$. The local ring $\hat {\Cal O}_{x,\Cal P_i}$ is a complete 
discrete valuation ring with uniformiser $\pi$.
We refer to $\{\Cal X_{i}\}_{i=1}^n$ as the set of {\it boundaries} of the formal germ $\Cal F_x$. 
We have a canonical morphism $\Cal X_{i}\to \Cal F_x$ of formal schemes, $1\le i\le n$.

Let $Z$ be a finite set of closed points of $X$ and
$U\subset X$ a formal sub-scheme of $X$ whose special fibre is $U_k\defeq X_k\setminus Z$.

\definition
{Definition 1.2.2}
We use the notations above. A ($G$-)covering patching data for the pair $(X,Z)$
consists of the following.

\noindent
(i)\  A finite (Galois) covering $V\to U$ of formal schemes (with Galois group $G$).

\noindent
(ii)\ For each point $x\in Z$, a finite (Galois) covering $\Cal Y_x\to \Cal F_x$ of formal schemes (with Galois group $G$).

The above data (i) and (ii) must satisfy the following compatibility condition.

\noindent
(iii)\ If $\{\Cal X_{i}\}_{i=1}^n$ are the boundaries of the formal germ at the point $x$, then for $1\le i\le n$
is given a ($G$-equivariant) $\Cal X_{i}$-isomorphism 
$$\Cal Y_x\times _{\Cal F_x}\Cal X_{i}\isom V\times _{U} \Cal X_{i}.$$
Property (iii) should hold for each $x\in Z$. (Note that there are natural morphisms $\Cal X_{i}\to U$ of formal schemes, $1\le i\le n$.)
\enddefinition

The following is the main patching result that we will use in this paper (cf. [Pries], Theorem 3.4, [Harbater], Theorem 3.2.8).

\proclaim {Proposition 1.2.3} We use the notations above. Given a ($G$-)covering patching data as in Definition 1.2.2 
there exists a unique, up to isomorphism, (Galois) covering $Y\to X$ of formal schemes (with Galois group $G$) which induces the 
above ($G$-)covering
in Definition 1.2.2(i) when restricted to $U$, and induces the above ($G$-)covering in Definition 1.2.2(ii) 
when pulled-back to $\Cal F_x$  for each point $x\in Z$.
\endproclaim

\subhead
{1.2.4}  
\endsubhead
With the same notations as above, let $x\in X$ be a closed point and $\widetilde X_k$ the normalisation of $X_k$. There is a one-to-one correspondence
between the set of points of $\widetilde X_k$ above $x$ and the set of boundaries of the formal germ of $X$ at the point $x$. Let $x_i$ be the point
of $\widetilde X_k$ above $x$ which corresponds to the boundary $\Cal X_{i}$, $1\le i\le n$. 
Then the completion of $\widetilde X_k$ at $x_i$ is isomorphic to the spectrum of a ring of formal power series $k[[t_i]]$ over $k$
where $t_i$ is a local parameter at $x_i$.  
The complete local ring $\hat {\Cal O}_{x,\Cal P_i}$ is a discrete valuation ring with uniformiser $\pi$ and residue field isomorphic
to $k((t_i))$. Fix an isomorphism $ k((t_i))\isom \hat {\Cal O}_{x,\Cal P_i}/\pi$.
Let $T_i\in \hat {\Cal O}_{x,\Cal P_i}$ be an element which lifts (the image in $\hat {\Cal O}_{x,\Cal P_i}/\pi$ under the above isomorphism of)
$t_i$; we shall refer to such an element $T_i$ as a parameter of 
$\hat {\Cal O}_{x,\Cal P_i}$, or of the boundary $\Cal X_{i}$.
Then there exists an isomorphism $R[[T_i]]\{T_i^{-1}\}\isom    \hat {\Cal O}_{x,\Cal P_i}$, where  
$$R[[T]]\{T^{-1}\} \defeq \Big\{\sum _{i=-\infty}^{\infty}a_iT^i,\ \underset 
{i\to -\infty} \to  \lim \vert a_i \vert=0 \Big\}$$ 
and $\vert \ \vert$ is a normalised absolute value of $R$ (cf. [Bourbaki], $\S2$, 5).

\subhead
{1.3}
\endsubhead
Let $X$ be a normal and geometrically connected flat $R$-scheme (resp. $R$-formal affine scheme) 
whose special fibre is equidimensional of dimension 1, $F\subset X_K$ a finite set of closed points,
and $\eta$ a geometric point of $X$ above its generic point. 
Then $\eta$ determines an algebraic closure $\overline K$ of $K$
and we have an exact sequence of arithmetic fundamental groups 
$$1\to \pi_1(X_K\setminus F,\eta)^{\geo}\to \pi_1(X_K\setminus F,\eta)\to \Gal (\overline K/K)\to 1,$$ 
where $\pi_1(X_K\setminus F,\eta)^{\geo}\defeq \Ker \left(\pi_1(X_K\setminus F,\eta)\twoheadrightarrow \Gal (\overline K/K)\right)$ is the geometric fundamental group of 
$X_K$ with generic point $\eta$. 
(In case $X=\Spf A$ is formal affine we define $\pi_1(X_K,\eta)\defeq \pi_1(\Spec A_K,\eta)$ and similarly we define $\pi_1(X_K\setminus F,\eta)$, cf. [Sa\"\i di], 2.1.)

\definition
{Definition 1.3.1} Let $S,T\subset X_K$ be (possibly empty) finite sets of closed points (which we also view as reduced closed sub-schemes of $X_K$). Assume that the special fibre $X_k$ of $X$ is reduced. Let $I\defeq I_{X_k,T}\subset \pi_1(X_K\setminus (S\cup T),\eta)^{\geo}$ be the subgroup normally generated by the inertia subgroups above the generic points of $X_k$
and the pro-$p$ Sylow subgroups of the inertia groups above all points in $T$. We define 
$$\widehat \pi_1(X_K\setminus (S\cup T);T,\eta)^{\geo}\defeq \pi_1(X_K\setminus (S\cup T),\eta)^{\geo}/I$$ 
and refer to it as the geometric \'etale fundamental group of $X_K\setminus (S\cup T)$; with base point $\eta$, {\it generically \'etale above $X_k$} and {\it tamely ramified} above $T$. In case $T=\emptyset$ and $U_K\defeq X_K\setminus S$ we simply write $\widehat \pi_1(U_K,\eta)^{\geo}\defeq \widehat \pi_1(X_K\setminus S;\emptyset,\eta)^{\geo}$. 
\enddefinition

Note that the definition of $\widehat \pi_1(X_K\setminus (S\cup T);T,\eta)^{\geo}$ depends on the model $X$ of $X_K$ (the model $X$ of $X_K$ will be fixed in later discussions in this paper). 
The profinite group $\widehat \pi_1(X_K\setminus (S\cup T);T,\eta)^{\geo}$ classifies finite covers 
$f:Y_L\to X_L\defeq X\times_KL$ where $L/K$ is a finite extension with valuation ring $R_L$, which are \'etale above $X_L\setminus (S\cup T)_L$ (here $(S\cup T)_L\defeq (S\cup 
T)\times _KL$) and possibly ramified with ramification indices prime-to-$p$
above the points in $T_L\defeq T\times_KL$, $Y_L$ is geometrically connected, and such that $f$ extends after possibly a finite extension of $L$
to a finite cover $\tilde f:Y\to X_{R_L}\defeq X\times _RR_L$ with $Y$ normal and $\tilde f$ is \'etale above the generic points of $X_k$.
Note that if $X$ is a smooth 
$R$-formal affine scheme as above which is an $R$-formal curve then $\widehat \pi_1(X_K,\eta)^{\geo}$ is isomorphic to the geometric \'etale fundamental group of the affine scheme $X_k$ 
as follows from the theorems of liftings of \'etale coverings (cf. [Grothendieck], Expos\'e I, Corollaire 8.4) and the theorem of purity of Zarizski-Nagata 
(cf. loc. cit. Expos\'e X, Th\'eor\`eme de puret\'e 3.1).
Note also that $\widehat \pi_1(X_K\setminus (S\cup T);T,\eta)^{\geo,p'}=\pi_1(X_K\setminus (S\cup T),\eta)^{\geo,p'}$, as follows easily from Abhyankar's lemma (cf. loc. cit. Expos\'e X, Lemme 3.6).

\subhead
\S2. Geometric Galois groups of formal boundaries of formal germs of $p$-adic formal curves
\endsubhead
In this section we investigate the structure of  a certain quotient of the geometric Galois group of a formal boundary of a formal germ of a formal $R$-curve. The results in this section will be used in 
$\S3$.

Let $D\defeq \Spf R\langle \frac {1}{T} \rangle$ 
be the formal standard $R$-closed unit disc and 
$D_K\defeq \Sp K\langle \frac {1}{T} \rangle$
its generic fibre which is the standard rigid $K$-closed unit disc centred at $\infty$.
Write $\Cal X=\Spf R[[T]]\{T^{-1}\}$ and $\Cal X_K\defeq \Spec (R[[T]]\{T^{-1}\}\otimes _RK)$. 
We have natural morphisms $\Cal X\to D$ of formal $R$-schemes, and $\Cal X_K\to D_K$ of $K$-schemes (cf. Notations).
Let $\eta$ be a geometric point of $\Cal X_K$ with value in its generic point which determines a generic point of 
$D_K$; which we denote also $\eta$, as well as algebraic closures $\overline K$ of $K$, $\overline k$ of $k$, and $\overline M$ of $M\defeq \Fr (R[[T]]\{T^{-1}\})$.
We have an exact sequence of Galois groups
$$1\to \Gal (\overline M/\overline K.M)\to \Gal (\overline M/M)\to 
\Gal (\overline K/K)\to 1.$$

Let $I\defeq I_{(\Cal X_k)}\subset \Gal (\overline M/\overline K.M)$ be the subgroup normally generated by the inertia subgroups above the generic point of $\Cal X_k$. Write 
$\Delta \defeq \Gal (\overline M/\overline K.M)/I$ and $\Gamma\defeq \Delta ^{p'}$. 
We have an exact sequence
$$1\to P\to \Delta \to \Gamma\to 1,$$
where $P\defeq \Ker (\Delta \twoheadrightarrow \Gamma)$.

\proclaim {Lemma 2.1} With the notations above, $P$ is the unique pro-$p$ Sylow subgroup of $\Delta$, $P$ is pro-$p$ free, and $\Gamma$ 
is canonically isomorphic to the maximal prime-to-$p$ quotient $\hat \Bbb Z (1)^{p'}$
of the Tate twist $\hat \Bbb Z(1)$.
\endproclaim

\demo{Proof} Indeed, it follows from the various Definitions that $\Delta$ is isomorphic to the absolute Galois group of $\overline k((t))$ which is known to be an extension
of $\hat \Bbb Z (1)^{p'}$ by a free pro-$p$ group.
\qed
\enddemo

\proclaim{Lemma 2.2}
Assume $p>0$. Then the pro-$p$ group $\widehat \pi_1(D_K,\eta)^{\geo,p}$ is free.
\endproclaim

\demo{Proof} Indeed, it follows from the various Definitions that $\widehat \pi_1(D_K,\eta)^{\geo,p}$ is isomorphic to the maximal pro-$p$ quotient of the geometric fundamental group of 
$D_k=\Bbb A^1_k$ which is pro-$p$ free (cf. [Serre1], Proposition 1).
\qed
\enddemo

\proclaim {Proposition 2.3} 
Assume $p>0$. Then the homomorphism  $\Delta \to \widehat \pi_1(D_{K},\eta)^{\geo}$ induced by the morphism
$\Cal X_K\to D_K$ induces a homomorphism 
$\Delta^{p}\to \widehat \pi_1(D_{K},\eta)^{\geo,p}$
which makes $\Delta ^{p}$ into a direct factor of $\widehat \pi_1(D_{K},\eta)^{\geo,p}$. Moreover,
$\Delta ^{p}$ is a free pro-$p$ group of infinite rank.
\endproclaim

\demo {Proof} We show that the map $\psi:H^1(\widehat \pi_1(D_{K},\eta)^{\geo},\Bbb Z/p\Bbb Z)
\to H^1(\Delta,\Bbb Z/p\Bbb Z)$ induced by the homomorphism  $\Delta \to \widehat \pi_1(D_{K},\eta)^{\geo}$ on cohomology is surjective 
(cf. Proposition 1.1.3). 
Let $\tilde f:\Delta\twoheadrightarrow \Bbb Z/p\Bbb Z$ be a surjective homomorphism and
$f:\Cal Y\to \Cal X$ the corresponding Galois cover (which we can assume, without loss of generality, defined over $K$) with $\Cal Y$ normal, geometrically connected, and $f$ is \'etale above the generic point of $\Cal X_k$ (hence $f$ is \'etale above $\Cal X$). Thus,
$f_k:\Cal Y_k\to \Cal X_k=\Spec k((t))$ is an \'etale $\Bbb Z/p\Bbb Z$-torsor. By Artin-Schreier theory the torsor $f_k$ can be approximated by a Galois cover 
$g_k:Y_k\to \Bbb P^1_k$ of degree $p$ which is \'etale outside the point $t=0$ and whose completion above this point is isomorphic to $f_k$.
The \'etale $\Bbb Z/p\Bbb Z$-torsor $g_k^{-1}(\Bbb A^1_k)\to \Bbb A^1_k\defeq \Bbb P^1_k\setminus \{t=0\}=\Spec k[\frac {1}{t}]$ lifts (uniquely up to isomorphism) to an \'etale $\Bbb Z/p\Bbb Z$-torsor $g:Z_K\to D_K$ by the 
theorems of liftings of \'etale covers
(cf. [Grothendieck], Expos\'e I, Corollaire 8.4) which gives rise to a class in  $H^1(\widehat \pi_1(D_{K},\eta)^{\geo},\Bbb Z/p\Bbb Z)$ that is easily verified to map to the class of $f$ in
$H^1(\Delta,\Bbb Z/p\Bbb Z)$. 
Moreover, $\Delta^p$ has infinite rank as it is isomorphic to the maximal pro-$p$ quotient of the absolute Galois group of $\overline k((t))$ which is known to be free of infinite rank.
\qed
\enddemo

Write $\widetilde \Gamma \defeq \widehat \pi_1(D_{K}\setminus \{\infty\}; \{\infty\},\eta)^{\geo,p'}=\pi_1(D_{K}\setminus\{\infty\},\eta)^{\geo,p'}$ (cf. 1.3) 
for the maximal prime-to-$p$ quotient of $\widehat \pi_1(D_{K}\setminus \{\infty\};\{\infty\},\eta)^{\geo}$.

\proclaim {Lemma 2.4} The morphism $\Cal X_K\to D_K$ induces a canonical homomorphism 
$\Gamma \to \widetilde \Gamma$ which is an isomorphism. In particular, $\widetilde \Gamma$ is (canonically) isomorphic to $\hat {\Bbb Z}(1)^{p'}$.
\endproclaim

\demo{Proof} Follows easily from the fact that a Galois covering $Y_K\to D_K$ of order prime-to-$p$ with $Y_K$
geometrically connected, ramified only above $\infty$ is, possibly after a finite extension of $K$ and for a suitable choice of the parameter $T$ of $D_K$,
generically a $\mu_n$-torsor given generically by the equation $S^n=T$
for some positive integer $n$ prime-to-$p$. 
\qed
\enddemo

Consider the following exact sequence
$$1\to {\Cal H}\to  \widehat \pi_1(D_{K}\setminus \{\infty\};\{\infty\},\eta)^{\geo}  \to \widetilde \Gamma \to 1,$$
where $\Cal H\defeq \Ker \lgroup \widehat \pi_1(D_{K}\setminus \{\infty\}; \{\infty\},\eta)^{\geo}  \twoheadrightarrow \widetilde \Gamma\rgroup$.
Further, let $\widetilde P\defeq {\Cal H}^p$ 
be the maximal pro-$p$ quotient of
${\Cal H}$. By pushing out the above sequence by the 
(characteristic) quotient ${\Cal H}\twoheadrightarrow \widetilde P$ we obtain an exact sequence
 $$1\to \widetilde {P}\to \widetilde \Delta \to \widetilde \Gamma \to 1.$$

\proclaim {Proposition 2.5} The morphism $\Cal X_K\to D_K$ induces a commutative diagram of exact sequences
$$
\CD
 1@>>> P @>>> \Delta @>>> \Gamma @>>> 1 \\
@.      @VVV     @VVV   @VVV \\
1@>>> \widetilde {P} @>>> \widetilde \Delta @>>> \widetilde \Gamma @>>> 1 \\
\endCD
$$
where the right vertical homomorphism $\Gamma\to \widetilde \Gamma$ is an isomorphism (cf. Lemma 2.4) 
and the middle vertical homomorphism $\Delta\to \widetilde \Delta$ makes $\Delta$ into a semi-direct factor of $\widetilde \Delta$ (cf. Lemma 1.1.5). 
\endproclaim

\demo{Proof} Let $\Delta \twoheadrightarrow G$ be a finite quotient which sits in an exact sequence $1\to Q\to G\to \Gamma_n\to 1$ where $\Gamma_n$
is the unique quotient of $\Gamma$ of cardinality $n$; for some integer $n$ 
prime-to-$p$, with $Q$ a $p$-group (cf. Lemma 2.1). We will show there exists a surjective homomorphism
$\widetilde \Delta\twoheadrightarrow G$ whose composition with $\Delta \to \widetilde \Delta$ is the above homomorphism.
We can assume, without loss of generality, that the corresponding Galois covering $f:\Cal Y\to \Cal X$ with group $G$ is defined over $K$, $\Cal Y$ is normal and connected, and $f$ is \'etale. 
This covering factorises as $\Cal Y\to \Cal X'\to \Cal X$ where  $\Cal X'\to \Cal X$ is Galois with group 
$\Gamma _n\isom \mu_n$ and $\Cal Y\to \Cal X'$ is Galois with group $Q$.  
After possibly a finite extension of $K$ the $\mu_n$-torsor $\Cal X'\to \Cal X$ extends to a generically $\mu_n$-torsor $D'\to D$ defined generically by an equation $S^n=T$, for a suitable choice of the parameter $T$ of $D$, which is (totally) ramified only above $\infty$, $D'=\Spf R\langle \frac{1}{S} \rangle$ is a closed formal unit disc centred at the unique point; which we denote 
also $\infty$, above $\infty \in D$ and $\Cal X'= \Spf R[[S]]\{S^{-1}\}$ (cf. Lemma 2.4 and the isomorphism $\Gamma\isom \widetilde \Gamma$ therein).

For the rest of the proof we assume $p>0$. By Proposition 2.3, applied to $\Cal X'\to D'$, there exists (after possibly a finite extension of $K$) an \'etale Galois covering 
$Y\to D'$ with group $Q$, $Y$ is normal and geometrically connected, and such that we have a commutative diagram of cartesian squares
$$
\CD
\Cal Y@>>> \Cal X ' @>>> \Cal X\\
@VVV @VVV   @VVV\\
Y @>>> D' @>>> D .\\
\endCD
$$
Next, we borrow some ideas from [Garuti] (preuve du Th\'eor\`eme 2.13). We claim one can choose the above (geometric) covering $Y\to D'$ 
such that the finite composite covering $Y\to D$ is Galois with group $G$. Indeed, consider the quotient
$\Delta \twoheadrightarrow \Delta _{\Cal X'}$ (resp. $\widetilde \Delta\twoheadrightarrow \widetilde \Delta_{D'}$) of $\Delta$ (resp. $\widetilde \Delta$)
which sits in the following exact sequence $1\to P_{\Cal X'}\to \Delta_{\Cal X'} \to \Gamma_n \to 1$ where $P_{\Cal X'}\defeq \widehat
\pi_1(\Cal X',\eta)^{\geo,p}$ (resp. $1\to \widetilde P_{D'}\to \widetilde \Delta_{D'} \to \widetilde \Gamma_n \to 1$ where $\widetilde P_{D'}\defeq \widehat
\pi_1(D',\eta)^{\geo,p}$). We have a commutative diagram of exact sequences

$$
\CD
1 @>>> P_{\Cal X'} @>>> \Delta_{\Cal X'} @>>> \Gamma_n @>>>1\\
@. @VVV @VVV @VVV\\
1 @>>> \widetilde P_{D'} @>>> \widetilde \Delta_{D'} @>>> \widetilde \Gamma_n @>>>1\\
\endCD
$$
where the right vertical map is an isomorphism (cf. Lemma 2.4). 
The choice of a splitting of the upper sequence in the above diagram (which splits since $P_{\Cal X'}$ is pro-$p$ and $\Gamma _n$ is cyclic (pro-)prime-to-$p$) induces
an action of $\Gamma_n$ on $\widetilde P_{D'}$ and $P_{\Cal X'}$ is a direct factor of $\widetilde P_{D'}$ (cf. Proposition 2.3) which is stable by this action of $\Gamma_n$. Further, $P_{\Cal X'}$ has a supplement $E$ in $\widetilde P_{D'}$ which is invariant under the action of $\Gamma_n$ by [Garuti], Corollaire 1.11. The existence of this supplement $E$ implies that one can choose 
$Y\to D'$ as above such that the finite composite covering $Y\to D$ is Galois with group $G$.
More precisely, if the Galois covering $\Cal Y\to \Cal X'$ corresponds to the surjective homomorphism $\rho:P_{\Cal X'}\twoheadrightarrow Q$ 
(which is stable by $\Gamma_n$ since $\Cal Y\to \Cal X$ is Galois) then we consider the Galois covering   $Y\to D'$ corresponding to the surjective homomorphism $\widetilde P_{D'}=P_{\Cal X'} \star E\twoheadrightarrow  Q$ which is induced by $\rho$ and the trivial homomorphism $E\to Q$, which is stable by $\Gamma_n$. 

The above construction can be performed in a functorial way with respect to the various finite quotients of $\Delta$. More precisely,
let $\{\phi_j:\Delta \twoheadrightarrow G_j\}_{j\in J}$ be a cofinal system of finite quotients of $\Delta$ where $G_j$
sits in an exact sequence $1\to Q_j\to G_j\to \Gamma_{n_j}\to 1$, for some integer $n_j$ 
prime-to-$p$, and $Q_j$ a $p$-group. Assume we have a factorisation $\Delta \twoheadrightarrow G_{j'}\twoheadrightarrow G_j$ for $j',j\in J$. Thus, $n_j$ divides $n_{j'}$, and we can assume without loss of generality (after replacing the group extension $G_j$ by its pull-back via $\Gamma _{n_j'}\twoheadrightarrow \Gamma _{n_j}$)
that $n\defeq n_j=n_{j'}$. With the above notations we then have surjective homomorphisms $\rho_{j'}:P_{\Cal X'}\twoheadrightarrow Q_{j'}$, 
$\rho_j:P_{\Cal X'}\twoheadrightarrow Q_j$ 
(which are stable by $\Gamma_n$), and $\rho_{j}$ factorises through $\rho_{j'}$. Then we consider the Galois covering(s)   $Y_{j'}\to D'$ (resp. 
$Y_j\to D'$) corresponding to the surjective homomorphism(s) $\psi_{j'}:\widetilde P_{D'}=P_{\Cal X'} \star E\twoheadrightarrow  Q$ 
(resp. $\psi_j:\widetilde P_{D'}=P_{\Cal X'} \star E\twoheadrightarrow  Q$) 
which are induced by $\rho_{j'}$ (resp. $\rho_j$) and the trivial homomorphism $E\to Q$, which are stable by $\Gamma_n$ and $\psi_{j}$ factorises through
$\psi_{j'}$. We deduce from this construction the existence of a surjective continuous homomorphism $\widetilde \Delta\twoheadrightarrow  \Delta$ which is a right inverse to the natural homomorphism $\Delta \to \widetilde \Delta$ (cf. Lemma 1.1.5).
\qed
\enddemo

\subhead
\S3 Geometric fundamental groups of formal fibres of $p$-adic curves
\endsubhead
In this section we investigate the structure of $\widehat \pi_1$ of a formal fibre of a $K$-curve.
Let $X$ be a proper, normal, formal $R$-curve whose special fibre $X_k$ is reduced and consists of $n\ge 1$ distinct
irreducible components $\{P_{i}\}_{i=1}^n$ which intersect at a (closed) point $x\in X_k(k)$, and $x$ is the unique singular point of $X_k$. 
Write $\widetilde P_i\to P_i$ for the morphism of normalisation. We assume $\widetilde P_i=\Bbb P^1_k$ is a projective line, the morphism $\widetilde P_i\to P_i$ is a homeomorphism, and if $x_i$ is the (unique) pre-image of $x$ in $\widetilde P_i$ then $x_i\in \widetilde P_i(k)$ is the zero point of $\widetilde P_i$. 
In particular, the configuration of the irreducible components of $X_k$ is tree-like. 
The formal curve $X$ has a formal covering $X=B\cup D_1\cup \ldots \cup D_n$ where
$B\subset X$ is a formal sub-scheme with special fibre $B_k=X_k\setminus \{\infty_i\}_{i=1}^n$, $\infty_i$ is the image in $P_i$ of the infinity point of 
$\widetilde P_i$,
$D_i=\Spf R\langle \frac {1}{T_i} \rangle$ is an $R$-formal closed unit disc with special fibre
$D_{i,k}=P_i\setminus \{x\}$ and generic fibre $D_{i,K}=\Sp K\langle \frac {1}{T_i} \rangle$ which is a closed unit $K$-rigid disc centred 
at the point $\infty_i\in D_{i,K}(K)$ which specialises into the infinity point $\infty_i\in P_i$, $1\le i\le n$.
Write $\Cal F\defeq {\Cal F}_x=\Spf \hat {\Cal O}_{X,x}$
for the formal germ of $X$ at $x$ and 
$\Cal F_K\defeq \Cal F_{x,K}=\Spec (\hat \Cal O_{X,x}\otimes _RK)$ for the formal fibre of $X_K$ at $x$ (cf. 1.2.1).
For $1\le i\le n$, let $\Cal X_i$ be the formal boundary of $\Cal F$ corresponding to the point $x_i$ above. The completion of the normalisation $X_k^{\nor}$ 
of $X_k$ at $x_i$ is isomorphic to the spectrum of a ring of formal power series $k[[t_i]]$ in one variable over $k$, and
$\Cal X_i\isom \Spf R[[T_i]]\{T_i^{-1}\}$ (cf. 1.2.4).

\definition {Remark 3.1}
Let $\widetilde Y$ be a proper and normal formal $R$-curve with $\widetilde Y_k$ reduced and $y\in \widetilde Y(k)$ a closed point. 
Write $\Cal G\defeq {\Cal G}_y=\Spf \hat {\Cal O}_{\widetilde Y,y}$
for the formal germ of $\widetilde Y$ at $y$ and 
$\Cal G_K\defeq \Cal G_{y,K}=\Spec (\hat \Cal O_{\widetilde Y,y}\otimes _RK)$ for the formal fibre of  
$\widetilde Y_K$ at $y$ (cf. 1.2.1).
Let $\{\Cal Y_i\}_{i=1}^n$ be the set of formal boundaries of $\Cal G$, and $y_i\in (\widetilde Y_k)^{\nor}(k)$ the point
of the normalisation $(\widetilde Y_k)^{\nor}$ of $\widetilde Y_k$ above $y$ which corresponds to the boundary $\Cal Y_{i}$, $1\le i\le n$.
The completion of $(\widetilde Y_k)^{\nor}$ at $y_i$ is isomorphic to the spectrum of a ring of formal power series $k[[s_i]]$ in one variable over $k$ and
$\Cal Y_i\isom \Spf R[[S_i]]\{S_i^{-1}\}$ (cf. 1.2.4). 
One can construct a compactification of $\Cal G$ (as in the above discussion where $\Cal G=\Cal F$)
which is a formal and proper $R$-curve $Y\defeq Y_y$ obtained by patching an $R$-formal closed unit disc 
$Y_i=\Spf R\langle \frac {1}{S_i} \rangle$ with $\Cal G$ along
the boundary $\Cal Y_i$, for $1\le i\le n$. 
The resulting formal $R$-curve $Y$ has a special fibre $Y_k$ consisting of $n$ distinct reduced
irreducible components 
$\{Q_{i}\}_{i=1}^n$ which intersect at the (closed) point $y$, and $y$ is the unique singular point of $Y_k$. Moreover, if we write $\widetilde Q_i\to Q_i$ for the morphism of normalisation then $\widetilde Q_i=\Bbb P^1_k$ is a projective line and the morphism $\widetilde Q_i\to Q_i$ is a homeomorphism. 
By construction the formal germ (resp. formal fibre) of $Y$ (resp. of $Y_K$) at the closed point $y$ is isomorphic to $\Cal G$ (resp. $\Cal G_{K}$).
(cf. [Bosch-L\"utkebohmert], Definition 4.4, for a rigid analytic construction of the generic fibre $Y_K^{\rig}\defeq Y_K$ of the above compactification $Y$ endowed with 
a formal covering corresponding to the above formal model $Y$ of $Y_K$, 
as well as [Bosch], Theorem 5.8, for the invariance of the formal germ at $y$ under this construction.) 
\enddefinition

Let $\eta$ be a geometric point of $\Cal F_{K}$ with value in its generic point which induces a geometric point $\eta$ of $X_K$ via the natural (scheme theoretic) 
morphism $\Cal F_{K}\to X_K$ (cf. Notations) and determines an algebraic closure $\overline K$ of $K$. 
Let $S\subset \Cal F_{K}$ be a (possibly empty) finite set of closed points. We have an exact sequence of arithmetic fundamental groups 
$$1\to \pi_1(\Cal F_{K}\setminus S,\eta)^{\geo}\to \pi_1(\Cal F_{K}\setminus S,\eta)\to \Gal (\overline K/K)\to 1,$$
where $\pi_1(\Cal F_{K}\setminus S,\eta)^{\geo}\defeq \Ker \lgroup \pi_1(\Cal F_{K}\setminus S,\eta)\twoheadrightarrow \Gal (\overline K/K)\rgroup$.
Write $\widehat \pi_1(\Cal F_{K}\setminus S,\eta)^{\geo}$ for the quotient of $\pi_1(\Cal F_{K}\setminus S,\eta)^{\geo}$ defined in 1.3.1.
Thus,  $\widehat \pi_1(\Cal F_{K}\setminus S,\eta)^{\geo}$
classifies finite (geometric) coverings $\Cal Y\to \Cal F$ (which we assume without loss of generality defined over $K$) with $\Cal Y$ normal and geometrically connected, 
which are \'etale above $\Cal F_K\setminus S$ and above the generic points of $\Cal F_k$. 
Write $U_{K}\defeq X_K\setminus (S\cup \{\infty_i\}_{i=1}^n)$ which is an affine curve and 
$\widehat \pi_1(U_K; \{\infty_i\}_{i=1}^n,\eta)^{\geo}\defeq \widehat \pi_1(X_K\setminus (S\cup\{\infty_i\}_{i=1}^n); \{\infty_i\}_{i=1}^n,\eta)^{\geo}$ 
for the quotient of  $\pi_1(U_K,\eta)^{\geo}$ defined in 1.3.1.
 Thus, $\widehat \pi_1(U_K; \{\infty_i\}_{i=1}^n,\eta)^{\geo}$ classifies finite (geometric) coverings $Y\to X$ (which we assume without loss of generality defined over $K$) with $Y$ 
normal and geometrically connected, which are \'etale above $X_K\setminus (S\cup\{\infty_i\}_{i=1}^n)$, are possibly ramified above the points $\{\infty_i\}_{i=1}^n$ with 
ramification indices prime-to-$p$, and are \'etale above the generic points of $X_k$. One of our main results is the following.

\proclaim {Theorem 3.2} The (scheme) morphism $\Cal F_{K}\to X_K$ induces a continuous homomorphism $\widehat \pi_1(\Cal F_{K}\setminus S,\eta)^{\geo}
\to \widehat \pi_1(U_K; \{\infty_i\}_{i=1}^n,\eta)^{\geo}$ (resp. $\widehat \pi_1(\Cal F_{K}\setminus S,\eta)^{\geo,p}
\to \widehat \pi_1(U_K; \{\infty_i\}_{i=1}^n,\eta)^{\geo,p}$)
which makes $\widehat \pi_1(\Cal F_{K}\setminus S,\eta)^{\geo}$ (resp. $\widehat \pi_1(\Cal F_{K}\setminus S,\eta)^{\geo,p}$) into a semi-direct factor of 
$\widehat \pi_1(U_K; \{\infty_i\}_{i=1}^n,\eta)^{\geo}$ (resp. $\widehat \pi_1(U_K; \{\infty_i\}_{i=1}^n,\eta)^{\geo,p}$). 
\endproclaim

\demo{Proof} We prove the first assertion by showing the criterion in Lemma 1.1.5 is satisfied. Let $\widehat \pi_1(\Cal F_K\setminus S,\eta)^{\geo}\twoheadrightarrow G$ be a finite quotient 
(which we can assume without loss of generality) corresponding to 
a finite Galois covering $f:\Cal Y\to \Cal F$ with group $G$, with $\Cal Y$ normal and geometrically connected,
which is \'etale above $\Cal F_K\setminus S$ and above the generic points of $\Cal F_k$. 
We will show the existence of a surjective homomorphism 
$\widehat \pi_1(U_K; \{\infty_i\}_{i=1}^n,\eta)^{\geo}\twoheadrightarrow G$ whose composite with  $\widehat \pi_1(\Cal F_K\setminus S,\eta)^{\geo}
\to \widehat \pi_1(U_K; \{\infty_i\}_{i=1}^n,\eta)^{\geo}$ is the above homomorphism.
For $1\le i\le n$, let $f_i:\Cal Y_i=\cup_{j=1}^{n_i} \Cal Y_{i,j}\to \Cal X_i$ be the pull-back 
of $f$ to $\Cal X_i$
via the natural morphism $\Cal X_i\to \Cal F$; $\{\Cal Y_{i,j}\}_{j=1}^{n_i}$ 
are the connected components of $\Cal Y_i$ and the morphism $f_{i,j}:\Cal Y_{i,j}\to \Cal X_i$ induced by $f_i$ is Galois with group $G_j$ a subgroup of $G$. Thus, $G_j$ is 
a quotient of $\widehat \pi_1(\Cal X_i,\eta_i)$ ($\eta_i$ is a suitable base point of $\Cal X_i$).
Fix $1\le j_0\le n_j$, then $f_i\isom \Ind _{G_{j_0}}^G f_{i,j_0}$ is an induced cover (cf. [Raynaud], 4.1).
By Proposition 2.5 there exists (after possibly a finite extension of $K$) a finite Galois covering $\tilde f_{i,j_0}:Y_{i,j_0}\to D_i$ with group $G_{j_0}$, where $Y_{i,j_0}$ is normal and geometrically connected, 
whose pull-back to $\Cal X_i$ via the natural morphism $\Cal X_i\to D_i$ is isomorphic to $f_{i,j_0}$. Further, the morphism 
$\tilde f_{i,j_0,}$ is ramified above $D_{i,K}$ possibly only above $\infty_i$ with ramification index prime-to-$p$, and $\tilde f_{i,j_0}$ is \'etale above the generic point of $D_{i,k}$.
Let $\tilde f_i : Y_i\defeq \Ind_{G_{j_0}}^G Y_{i,j_0}\to D_i$ be the induced cover (cf. loc. cit.), for $1\le i\le n$. By Proposition 1.2.3 
one can patch the covering $f$ with the coverings $\{\tilde f_i\}_{i=1}^n$ to construct a finite Galois covering  $\tilde f:Y\to X$ between formal $R$-curves 
with group $G$, $Y$ is normal and geometrically connected (since $\Cal Y_K$ is), which gives rise (via the formal GAGA functor) to a surjective homomorphism 
$\widehat \pi_1(U_K;\{\infty_i\}_{i=1}^n,\eta)^{\geo}\twoheadrightarrow G$
as required. 
Moreover, one verifies easily that the above construction can be performed in a functorial way with respect to the various quotients 
of $\widehat \pi_1(\Cal F_K\setminus S,\eta)^{\geo}$ (in the sense of lemma 1.1.5) using Proposition 2.5, so that one deduces the existence of a continuous homomorphism
$\widehat \pi_1(U_K;\{\infty_i\}_{i=1}^n,\eta)^{\geo}\to \widehat \pi_1(\Cal F_K\setminus S,\eta)^{\geo}$ which is right inverse to $\widehat \pi_1(\Cal F_K\setminus S,\eta)^{\geo}
\to \widehat \pi_1(U_K; \{\infty_i\}_{i=1}^n,\eta)^{\geo}$. The proof of the second assertion is entirely similar using similar arguments. 
\qed
\enddemo

\proclaim {Proposition 3.3} The (scheme) morphism $\Cal F_K\to X_K$ induces a continuous homomorphism $\widehat \pi_1(\Cal F_K\setminus S,\eta)^{\geo,p'}
\to \widehat \pi_1(U_K;\{\infty_i\}_{i=1}^n,\eta)^{\geo,p'}$ which makes $\widehat \pi_1(\Cal F_K\setminus S,\eta)^{\geo,p'}$ into a semi-direct factor of
$\widehat \pi_1(U_K;\{\infty_i\}_{i=1}^n,\eta)^{\geo,p'}$.
\endproclaim

\demo {Proof} The proof follows by using similar arguments to the ones used in the proof of Theorem 3.2. 
More precisely, with the notations in the proof of Theorem 3.2 the morphism $\Cal Y_{i,j}\to \Cal X_i$ in this case is Galois with group $\mu_n$, where $n$ is an integer prime-to-$p$, 
and extends (uniquely, possibly after a finite extension of $K$) 
to a cyclic Galois covering $Y_{i,j}\to D_i$ of degree $n$ ramified only above $\infty _i$ 
(cf. Lemma 2.1 and Lemma 2.4). 
\enddemo

In [Sa\"\i di1] we defined the {\it genus} $g_x$ of the closed point $x$ of $X$, whose definition depends only on the local (\'etale) structure of $X_k$ at $x$, and which equals the genus 
of the proper, connected, and smooth $K$-curve $X_K$ constructed above (cf. loc. cit. Lemma 3.3.1 and the discussion before it).
(The genus $g_x$ of $x$ is also called the genus of the formal fibre $\Cal F_K$.)

\proclaim {Theorem 3.4} Let $S(\overline K)=\{y_1,\ldots,y_m\}$ of cardinality $m\ge 0$. Then the continuous homomorphism 
$\pi_1(\Cal F_K\setminus S,\eta)^{\geo,p'}\to \pi_1(X_K\setminus (S\cup \{\infty_i\}_{i=1}^n),\eta)^{\geo,p'}$ (cf. Proposition 3.3) is an isomorphism. In particular,
$\pi_1(\Cal F_K\setminus S,\eta)^{\geo,p'}$ is (pro-)prime-to-$p$
free of rank $2g_x+n+m-1$ and can be generated by $2g_x+n+m$ generators 
$\{a_1,\ldots,a_g,b_1,\ldots,b_g,\sigma_1,\ldots,\sigma_n,\tau_1,\ldots,\tau_m\}$
subject to the unique relation $\prod_{i=1}^g[a_i,b_i]\prod _{j=1}^n\sigma _j\prod_{t=1}^m \tau _t=1$, where $\sigma_j$ (resp $\tau _t$) is a generator of inertia at $\infty _i$
(resp. $y_t$). 
\endproclaim

\demo{Proof} The homomorphism $\pi_1(\Cal F_K\setminus S,\eta)^{\geo,p'}\to \pi_1(X_K\setminus (\{\infty_i\}_{i=1}^n\cup S),\eta)^{\geo,p'}$ is injective as follows from Proposition 3.3 (note that $\widehat \pi_1=\pi_1$ in this case). We show it is surjective. To this end it suffices to show that given a finite Galois covering $f:Y\to X$ with group $G$ of cardinality prime-to-$p$, with $Y$ normal and geometrically connected, which is \'etale above $X_K\setminus (S\cup \{\infty_i\}_{i=1}^n)$, and $\tilde f:\Cal Y_K\to \Cal F_K$ its restriction to $\Cal F_K$, then $\Cal Y_K$ is geometrically connected. Equivalently, we need to show (possibly after passing to a finite extension of $K$) that $f^{-1}(x)$ consists of a single closed point.
(The set of connected components of $\Cal Y_K$ is in one-to-one correspondence with the set $f^{-1}(x)$.)
We can assume, without loss of generality, 
that $Y_k$ is reduced (cf. Lemme d'Abhyankar, [Grothendieck], Expos\'e X, Lemme 3.6). 
Let $y\in f^{-1}(x)$ and $D_y\subset G$ its decomposition group. Let $Y_i$ be an irreducible component of $Y_k$ above $P_{i}$ passing through $y$, $\widetilde Y_i\to Y_i$ the morphism of normalisation, and $\widetilde Y_i\to \widetilde P_i$ the natural morphism which is Galois
with group $D_{Y_i}\subset G$ the decomposition group of $Y_i$. The morphism $\widetilde Y_i\to \widetilde P_i$ is \'etale outside $\{x_i,\infty_i\}$ 
by Zariski's purity Theorem. 
Hence $D_{Y_i}=\mu_n$ is cyclic of order $n$, for some integer $n$ prime-to-$p$, and the above morphism $\widetilde Y_i\to \widetilde P_i$
is totally ramified above $\infty _i$ and $x_i$ as follows from the structure of $\pi_1(\Bbb P^1_{\bar k}\setminus \{0,\infty\})^{p'}$. In particular, $D_{Y_i}\subset D_y$. Moreover, 
$Y_k$ is regular outside $f^{-1}(x)$ (cf. [Raynaud], Lemma 6.3.2).
We can associate a graph $\Gamma$ to $Y_k$ whose vertices are the irreducible components of $Y_k$ and edges are the closed points of $Y_k$ above $x$, two vertices 
$Y_i$ and $Y_{i'}$ passing by a closed point $y$ above $x$ are linked by the edge $y$.
Assume that $f^{-1}(x)$ has cardinality $>1$ and let $\{y,y'\}\subseteq f^{-1}(x)$ be two distinct points. Then no irreducible component of $Y_k$ passes through both $y$ and $y'$
(cf. the above fact that $\widetilde Y_i\to \widetilde P_i$ is totally ramified above $x_i$). 
More precisely, if $Y_i$ is an irreducible component of $Y_k$ then $Y_i$ passes through a unique point $y$ of $Y_k$ above $x$.  From this (and the above facts) it follows easily that the connected components of $\Gamma$ 
are in one-to-one correspondence with the elements of $f^{-1}(x)$ and $\Gamma$
is disconnected which contradicts the fact that $Y_k$ is connected. Thus,  $f^{-1}(x)$ has cardinality $1$ necessarily as required.
The last assertion follows form the well-known structure of $\pi_1(X_K\setminus (S\cup \{\infty_i\}_{i=1}^n),\eta)^{\geo,p'}$ 
(cf. [Grothendieck], Expos\'e XIII, Corollaire 2.12). 
\qed
\enddemo

\definition {Examples 3.5} Suppose $K$ is of mixed characteristics with $\char (k)=p>0$.
Let $\Cal F=\Spf R[[T]]$ (resp. $\Cal F=\Spf R[[T_1,T_2]]/(T_1T_2-\pi^e)$) be the formal open unit disc (resp. formal open annulus of thickness $e\ge 1$) and 
$S=\{y_1,\ldots,y_m\} \subset \Cal F(K)$ a set of $m\ge 0$ distinct $K$-rational points (in the second case $e>1$ necessarily if $m\neq 0$). 
In this case $\Cal F$ has a compactification $X=\Bbb P^1_{R}$ the $R$-projective line with parameter $T$ and $\Cal F$ is the formal germ at $T=0$
(resp. a compactification $X$ which is a formal model of the projective line $\Bbb P^1_K$ consisting
of two formal closed unit discs $D_1$ and $D_2$ centred at $\infty _1$ and $\infty_2$; respectively, which are patched with $\Cal F$ along its two boundaries. 
The special fibre $X_k$ consists of two projective lines which intersect at the double point $x$ and $\Cal F$ is the formal germ at $x$).
The results of $\S3$ in this case read as follows. First, the homomorphism $\widehat \pi_1(\Cal F_K\setminus S,\eta)^{\geo}
\to \widehat \pi_1(\Bbb P^1_K\setminus (T\cup \{\infty\}); \{\infty\},\eta)^{\geo}$ (resp. $\widehat \pi_1(\Cal F_K\setminus S,\eta)^{\geo}
\to \widehat \pi_1(\Bbb P^1_K\setminus (T\cup \{\infty_1,\infty_2\}); \{\infty_1,\infty_2\},\eta)^{\geo}$) makes $\pi_1(\Cal F_K\setminus S,\eta)^{\geo}$ 
into a semi-direct factor of $\widehat \pi_1(\Bbb P^1_K\setminus (T\cup \{\infty\}); \{\infty\},\eta)^{\geo}$ (resp. 
$\widehat \pi_1(\Bbb P^1_K\setminus (T\cup \{\infty_1,\infty_2\}); \{\infty_1,\infty_2\},\eta)^{\geo}$) and the maximal prime-to-$p$ quotient $\pi_1(\Cal F_K\setminus S,\eta)^{\geo,p'}$
is free of rank $m$ (resp. $m+1$).
\enddefinition

$$\text{References.}$$



\noindent 
[Bosch] S. Bosch, Eine bemerkenswerte Eigenshaft des formellen Fasern affinoider R\"aume, Math. Ann. 229 (1977), 25--45.

\noindent
[Bosch-L\"utkebohmert] S. Bosch and W. L\"utkebohmert, Stable reduction and uniformisation of abelian varieties I, Math. Ann. 270 (1985), 349--379.


\noindent
[Bourbaki] N. Bourbaki, Alg\`ebre Commutative, Chapitre 9, Masson, 1983.


\noindent
[Garuti] M. Garuti, Prolongements de rev\^etements galoisiens en g\'eom\'etrie 
rigide, Compositio Mathematica, 104 (1996), no 3, 305--331.

\noindent
[Grothendieck] A. Grothendieck, Rev\^etements \'etales et groupe fondamental, Lecture 
Notes in Math. 224, Springer, Heidelberg, 1971.

\noindent
[Harbater] D. Harbater, Galois groups and fundamental groups, 313-424, Math. Sci. Res. Inst. Publ., 41, Cambridge Univ. Press, Cambridge, 2003.

\noindent
[Raynaud] M. Raynaud, Rev\^etements de la droite affine en caract\'eristique $p>0$
et conjecture d'Abhyankar, Invent. Math. 116 (1994), no 1-3, 425--462.

\noindent
[Ribes-Zalesskii] L. Ribes and P. Zalesskii, Profinite groups, Ergebnisse der Mathematik 
und ihrer Grenzgebiete. Folge 3. A series of Modern Survey in Mathematics 40. Springer-Verlag, Berlin 2000.

\noindent
[Pries] R. Pries, Construction of covers with formal and rigid geometry, in: J. -B. Bost, F. Loeser, M. Raynaud (Eds.), Courbes semi-stables et groupe fondamental en g\'eom\'etrie alg\'ebrique, Progr. Math., Vol. 187, 2000.

\noindent
[Sa\"\i di] M. Sa\"\i di, \'Etale fundamental groups of affinoid $p$-adic curves, Journal of algebraic geometry, 27 (2018), 727--749.

\noindent
[Sa\"\i di1] M. Sa\"\i di, Wild ramification and a vanishing cycles formula, 
J. Algebra 273 (2004), no. 1, 108--128.  



\noindent
[Serre] J-P. Serre, Cohomologie Galoisienne, Lecture Notes in Math., 5, Springer Verlag, 
Berlin, 1994.

\noindent
[Serre1] J-P. Serre, Construction de rev\^etements \'etale de la droite affine en caract\'eristique $p>0$, C. R. Acad. Sci. Paris 311 (1990), 341--346.


\bigskip

\noindent
Mohamed Sa\"\i di

\noindent
College of Engineering, Physics and Mathematical Sciences

\noindent
University of Exeter

\noindent
Harrison Building

\noindent
North Park Road

\noindent
EXETER EX4 4QF 

\noindent
United Kingdom

\noindent
M.Saidi\@exeter.ac.uk

\end
\enddocument